%% file: robust_erm.tex
\documentclass[12pt]{article}
\usepackage{fullpage,graphicx,psfrag,amsmath,amsfonts,verbatim}
\usepackage[small,bf]{caption}\usepackage{amssymb}
\usepackage{enumitem}
\usepackage{tikz}
\usepackage{listings}
\usepackage{authblk}
\usepackage{textcomp}
\usepackage{comment}

\graphicspath{{figures/}}

\definecolor{codegreen}{rgb}{0,0.6,0}
\definecolor{codegray}{rgb}{0.5,0.5,0.5}
\definecolor{codepurple}{rgb}{0.58,0,0.82}
\definecolor{backcolour}{rgb}{0.95,0.95,0.98}

\lstdefinestyle{mystyle}{
backgroundcolor=\color{backcolour},   
keywordstyle=\color{magenta},
numberstyle=\tiny\color{codegray},
stringstyle=\color{codepurple},
basicstyle=\ttfamily\small,
breakatwhitespace=false,         
breaklines=true,                 
captionpos=t,                    
keepspaces=true,                 
numbers=left,                    
numbersep=5pt,                  
showspaces=false,                
showstringspaces=false,
showtabs=false,                  
tabsize=2,
}

\definecolor{seagreen}{rgb}{0.18, 0.55, 0.34}
\definecolor{mediumviolet-red}{rgb}{0.78, 0.08, 0.52}
\definecolor{khaki}{rgb}{0.94, 0.9, 0.55}

\lstdefinelanguage{mypython}
{
keywords=[1]{from, import, as, assert, not, print, nonneg, PSD, for, in, return},
keywordstyle=[1]{\color{mediumviolet-red}},
keywords=[2]{surecr, torch, cp, lo, pl, cvxpy, dsp, Variable, LocalVariable,
sqrt, exp, saddle_inner, saddle_max, saddle_min, SaddlePointProblem, MinimizeMaximize,
numpy, np, Problem, Minimize, Maximize, is_dsp, value, solve, inner,
convex_variables, concave_variables, affine_variables, sum, multiply,
arange, norm1, norm2, norm_inf, abs, square, saddle_quad_form,
sum_squares, SuppFunc, transforms, suppfunc, transforms,
diagonal, outer, typing, Constraint, ndarray, ndim, power},
keywordstyle=[2]{\color{seagreen}},
keywords=[3]{append, range, NotImplementedError, List, Callable, len, str, def},
keywordstyle=[3]{\color{cyan}},
upquote=true,
showstringspaces=false,
basicstyle=\ttfamily,
columns=fullflexible,
keepspaces=true,
emphstyle={\color{seagreen}},
belowskip=1em,
aboveskip=1em,
morecomment=[l]{\#}
}

\usepackage{hyperref}
\usepackage{wasysym}
\hypersetup{
colorlinks   = true,    
urlcolor     = blue,    
linkcolor    = blue,    
citecolor    = red      
}

\usepackage[
backend=bibtex,
style=alphabetic,
sorting=ynt,
firstinits=true,
url=false, doi=false, eprint=false,
]{biblatex}
\usepackage{color}
\input defs.tex

\title{Specifying and Solving Robust Empirical Risk Minimization Problems Using CVXPY} 

\author[1]{Eric Luxenberg\footnote{Equal contribution.}}
\newcommand\CoAuthorMark{\footnotemark[\arabic{footnote}]} 
\author[2]{Dhruv Malik\protect\CoAuthorMark}
\author[2]{\\Yuanzhi Li}
\author[2]{Aarti Singh}
\author[1]{Stephen Boyd}
\affil[1]{Department of Electrical Engineering, Stanford University}
\affil[2]{Machine Learning Department, Carnegie Mellon University}

\begin{document}
\maketitle

\begin{abstract}
We consider robust empirical risk minimization (ERM), where model parameters are
chosen to minimize the worst-case empirical loss when each data point varies
over a given convex uncertainty set. In some simple cases, such problems can
be expressed in an analytical form. In general the problem can be made tractable
via dualization, which turns a min-max problem into a min-min problem.
Dualization requires expertise and is tedious and error-prone. 
We demonstrate how CVXPY can be used to automate this dualization procedure
in a user-friendly manner.
Our framework allows practitioners to specify and solve robust ERM problems with a
general class of convex losses, capturing many standard regression and
classification problems. Users can easily specify any complex uncertainty
set that is representable via disciplined convex programming (DCP) constraints.
\end{abstract}

\section{Robust empirical risk minimization}

Robust optimization is used in mathematical optimization,
statistics, and machine learning, to
handle problems where the data is uncertain.  In this note we
consider the robust empirical risk minimization (RERM) problem 
\BEQ\label{e-main}
\begin{array}{ll} \mbox{minimize} & \sum_{i=1}^n\sup_{x_i \in
\mathcal{X}_i}f(x_i^T\theta-y_i) \\
\mbox{subject to} & \theta \in \Theta,
\end{array}
\EEQ
with variable $\theta \in \reals^d$.
Here, $\Theta \subseteq \reals^d$ is closed and convex,
$\mathcal{X}_i \subset \reals^{d}$ is compact and convex for each
$i=1,\ldots,n$, $f:\reals\to\reals$ is convex and $\{ x_i, y_i \}_{i=1}^n$
is a dataset.
The objective is to
find $\theta \in \Theta$ that minimizes the worst-case value of
$\sum_{i=1}^nf(x_i^T\theta-y_i)$ over all possible $x_i$ in the given uncertainty
sets $\mathcal{X}_i$. Beyond convexity, we will assume that $f$ is
either non-increasing, or $f$ is non-decreasing on $\reals_+$ and a
function of the absolute value of its argument,
\ie, $f(z)=f(|z|)$.
\paragraph{Examples.}
Our assumptions capture a wide range of loss functions in both
regression and classification, including the following.
\begin{itemize}
\item \emph{Finite $p$-norm loss.} $f(z) = \vert z\vert^p$ for
$1 \leq p < \infty$. 
\item \emph {Huber loss.} $f(z) = \frac{1}{2}z^2$ for  $|z| \leq
\delta$, and $f(z)= \delta|z|-\frac{\delta^2}{2}$ for $|z| > \delta$, where 
$\delta>0$ is a parameter.
\item \emph{Hinge loss.} $f(z) = \max(0,1-z)$.
\item \emph{Logistic loss.} $f(z) = \log(1+\exp(-z))$.
\item \emph{Exponential loss.} $f(z) = \exp(-z)$.
\end{itemize}
Our formulation includes the case
of using hinge, logistic, or exponential loss for
binary classification, by solving \eqref{e-main} with the
transformed dataset $\{ y_i x_i, 0 \}$.

\subsection{Solving RERM problems}
\label{sec:solving_rerm_intro}
The problem \eqref{e-main} is convex, but not immediately tractable
because of the suprema appearing in the worst-case loss terms.
It can often be transformed to an explicit tractable form that does not
include suprema.

\paragraph{Analytical cases.} In some simple cases we can directly 
work out a tractable expression for the worst-case loss.
As a simple example, consider
$\mathcal X_i = \{x_i \mid \| x_i - \tilde x_i\|_2\leq \rho\}$,
where $\rho>0$.
When $f$ is non-increasing, the worst-case loss term is
\[
\sup_{x_i \in \mathcal{X}_i}f(x_i^T\theta-y_i) = 
f(\tilde x_i^T\theta -y_i - \rho\|\theta\|_2).
\]
When $f$ is non-decreasing on $\reals_+$ with $f(z)=f(|z|)$,
the worst-case loss term is
\[
\sup_{x_i \in \mathcal{X}_i}f(x_i^T\theta-y_i) = 
f(|\tilde x_i^T\theta -y_i| + \rho\|\theta\|_2).
\]
Both righthand sides are explicit convex expressions
that comply with the disciplined convex programming (DCP) rules.
This means they can be directly typed into domain specific languages (DSLs)
for convex optimization such as CVXPY \cite{diamond2016cvxpy}.

\paragraph{Dualization.} For more complex uncertainty sets the problem \eqref{e-main}
can still be transformed to a tractable form, using
\emph{dualization} of the suprema apprearing in the worst-case loss terms.
This dualization process converts the suprema in \eqref{e-main} to infima,
so that the problem can be solved by standard methods as a single minimization
problem.
Unfortunately, this dualization procedure is
cumbersome and error-prone. Many
practitioners are not well versed in this procedure, limiting its
use to experts. Moreover, a key step in this procedure
involves writing down a conic representation of
$\mathcal{X}_i$.  Such a calculation is antithetical to the
spirit of DSLs such as CVXPY, which were
introduced precisely to alleviate users of this burden.

\paragraph{Automatic dualization via CVXPY.}
In this note we show how
CVXPY can be used to conveniently solve \eqref{e-main}
with just a few lines of code,
even when the uncertainty sets $\mathcal{X}_i$ are complicated.
We also demonstrate how DSP, a recent DSL for disciplined
saddle programming \cite{schiele2023dsp} that is based
on CVXPY, can solve the RERM problem \eqref{e-main}
with the same ease and convenience.  In both approaches
no explicit dualization is needed, and
the code is short and naturally follows the math.
We demonstrate our approach with a synthetic 
regression example that, however,
uses real data, where the uncertainty sets are 
intervals intersected with a Euclidean ball.


\subsection{Previous and related work}
\paragraph{Robust optimization and saddle problems.}
Robust optimization is an approach that takes into account uncertainty,
variability or missing-ness of problem parameters 
\cite{ben2009robust}.
Saddle problems are robust optimization problems that include
the partial supremum or infimum of convex-concave saddle
functions. While \eqref{e-main} is not a priori a
saddle problem, we can solve it via DSP \cite{schiele2023dsp}, a recently
introduced DSL for saddle programming.

\paragraph{RERM.}
In machine learning and statistics, it is common to
learn a robust predictor or classifier by solving
\eqref{e-main} with
appropriate choices of $f, \mathcal{X}_i$
\cite{el1997robust, xu2008lasso,
xu2009svmrobust, bertsimas2011theory}.
When each $\mathcal{X}_i$ has benign structure,
then \eqref{e-main} admits
convenient reformulation for many choices of $f$ \cite{cvxbook}.
As an example, such reformulations have been
applied to learn linear regression functions 
when the feature matrix has missing data, and the features are known
to lie with high probability in an ellipsoid, so that \eqref{e-main} is easily
written as an SOCP \cite{shivaswamy2006socp, aghasi2022rigid}.
However, when $\mathcal{X}_i$ is not a simple set such as an
ellipsoid or box, then prior techniques reformulate
\eqref{e-main} by writing
$\mathcal{X}_i$ in conic form and then
dualizing \cite{ben2009robust}.

\section{Reformulating the RERM problem}
\label{s-main-results}
Throughout, our only requirement on the uncertainty sets
$\mathcal{X}_i$ is that each is a compact, convex set that can be
expressed via DCP constraints. This includes canonical scenarios,
such as when $\mathcal{X}_i$ is a polytope, or is a norm ball
centered at a nominal value.
But it also includes many complex uncertainty sets,
such as the intersection of a norm ball and a
polytope. We now reformulate \eqref{e-main}
in a manner that permits easy specification and solution
via CVXPY, under various monotonicity assumptions on $f$. 
Recall that
the support function of a non-empty closed convex set $C$ is given by
$\mathcal{S}_C(\theta)=\sup \{ x^T \theta : x \in C \}$, 
which is a fundamental object in convex analysis \cite{cvxbook}.

Introducing the epigraph variables $c \in
\reals^n$, the problem \eqref{e-main} is straightforwardly equivalent
to 
\BEQ\label{e-rerm-reform-1}
\begin{array}{ll} \mbox{minimize} & \sum_{i=1}^n c_i \\
\mbox{subject to} & \theta \in \Theta, \\
& \sup_{x_i \in
\mathcal{X}_i} f(x_i^T\theta - y_i) \leq c_i, \quad i = 1, \ldots, n. 
\end{array}
\EEQ
with variables $c \in \reals^n, \theta \in \reals^d$. We now
use the assumptions on $f$ to
rewrite the constraints $\sup_{x_i \in \mathcal{X}_i}
f(x_i^T\theta - y_i) \leq c_i$ in a tractable form.

\paragraph{Loss $f$ is non-increasing.} If $f$ is non-increasing, then
introducing an auxiliary variable $z_i$ shows that
\begin{eqnarray*}
\sup_{x_i \in \mathcal{X}_i} f( x_i^T\theta - y_i) \leq c_i
&\iff& f \left(\inf_{x_i \in \mathcal{X}_i} x_i^T\theta - y_i \right) \leq c_i \\
&\iff& \inf_{x_i \in \mathcal{X}_i} x_i^T\theta - y_i \geq z_i, \quad f(z_i) \leq c_i \\
&\iff& \sup_{x_i \in \mathcal{X}_i} -x_i^T\theta + y_i \leq -z_i, \quad f(z_i) \leq c_i.
\end{eqnarray*}
So, after eliminating the epigraph variable $c$ from
\eqref{e-rerm-reform-1}, we have shown
\eqref{e-main} is equivalent to
\BEQ\label{e-tractable-classification}
\begin{array}{llll} \mbox{minimize} & \sum_{i=1}^n f(z_i) \\
\mbox{subject to} & \theta \in \Theta,\\
& \mathcal{S}_{\mathcal{X}_i}(-\theta) + y_i \leq -z_i, \quad i = 1,
\ldots, n,
\end{array}
\EEQ
with variables $z \in \reals^n, \theta \in \reals^d$.
Typical classification losses, such as the hinge, logistic and exponential
losses, are non-increasing.

\paragraph{Loss $f$ is non-decreasing on $\reals_+$ and $f(a) = f(\vert a \vert)$.}
If $f$ is monotone on nonnegative arguments, and depends only on its
argument through the absolute value, then introducing the auxiliary
variable $z_i$ shows that
\begin{eqnarray*}
\sup_{x_i \in \mathcal{X}_i} f( x_i^T\theta - y_i) \leq c_i
&\iff& f \left( \sup_{x_i \in \mathcal{X}_i} \vert x_i^T\theta - y_i \vert \right) \leq c_i \\
&\iff& \sup_{x_i \in \mathcal{X}_i} \vert x_i^T\theta - y_i \vert \leq z_i, \quad f( z_i ) \leq c_i \\
&\iff& \sup_{x_i \in \mathcal{X}_i} x_i^T\theta - y_i \leq z_i,
\quad \sup_{x_i \in \mathcal{X}_i} -x_i^T\theta + y_i \leq z_i,
\quad f( z_i ) \leq c_i.
\end{eqnarray*}
So, after eliminating the epigraph variable $c$ from
\eqref{e-rerm-reform-1}, we have shown
\eqref{e-main} is equivalent to
\BEQ\label{e-tractable-regression}
\begin{array}{llll} \mbox{minimize} & \sum_{i=1}^n f(z_i) \\
\mbox{subject to} & \theta \in \Theta,\\
& \mathcal{S}_{\mathcal{X}_i}(\theta) - y_i \leq z_i, \quad i = 1, \ldots, n, \\
& \mathcal{S}_{\mathcal{X}_i}(-\theta) + y_i \leq z_i, \quad i = 1, \ldots,
n,
\end{array}
\EEQ
with variables $z \in \reals^n, \theta \in \reals^d$.
Typical regression losses, such as $p$-norm and Huber losses, satisfy
this requirement on $f$.

\paragraph{CVXPY code.}
\label{sec:results_cvxpy_code}
The robust constraints in
\eqref{e-tractable-classification} and \eqref{e-tractable-regression}
include suprema over $\mathcal{X}_i$ of bilinear forms involving
$x_i, \theta$. While this ostensibly requires dualization to handle,
the CVXPY transform \verb|SuppFunc| allows one to easily specify the
support function $\mathcal{S}_{C}(\theta)$ of a set $C$ created via DCP constraints.
Since this function is already implemented in CVXPY,
we can directly specify the robust constraints in
\eqref{e-tractable-classification} and \eqref{e-tractable-regression}
without additional reformulation or dualization.
As an example, we depict below
the CVXPY code that specifies and solves
\eqref{e-tractable-regression} with $f = \vert \cdot \vert^2$
and $\Theta = \reals^d$, which is a robust least squares
problem. For convenience, we assume $y$
has already been specified as \verb|y|.

\begin{lstlisting}[language=mypython]
import cvxpy as cp
from cvxpy.transforms.suppfunc import SuppFunc

theta, z = cp.Variable(d), cp.Variable(n)
constraints = []

for i in range(n):
    # Create variables for uncertainty set
    x = cp.Variable(d)
    
    # Construct uncertainty set containing x (filled in by user)
    local_constraints = []
    
    # Implement the support function of the uncertainty set
    G1 = SuppFunc(x, local_constraints)(theta)
    G2 = SuppFunc(x, local_constraints)(-theta)
    
    # Store robust constraints
    constraints.append(G1 - y[i] <= z[i])
    constraints.append(G2 + y[i] <= z[i])

obj = cp.Minimize(cp.sum_squares(z))
prob = cp.Problem(obj, constraints)
prob.solve()
\end{lstlisting}
To fully specify the problem, the user only needs to describe the uncertainty
set $\mathcal{X}_i$ for each $i = 1, \ldots, n$ in Line 12, in terms of
the \verb|x| instantiated in Line 9. This is done
exactly as one would typically do for any \verb|Variable| in CVXPY. If,
for example, $\mathcal{X}_i$ was the intersection of the Euclidean
unit ball, the non-negative orthant, and the set of vectors whose first
coordinate is 0.25, then replacing Line 12 with
the code block below is sufficient.
\begin{lstlisting}[language=mypython,firstnumber=12]
local_constraints = [x >= 0, cp.sum_squares(x) <= 1, x[0] == 0.25]
\end{lstlisting}
This manner of expressing $\mathcal{X}_i$ is thus natural, user-friendly
and directly follows the math.

Alternatively, one may recognize that the constraints in
\eqref{e-tractable-classification} and \eqref{e-tractable-regression}
include the partial suprema of a convex-concave saddle function.
Since DSP was designed to solve saddle problems, and 
a bilinear function is an atom in DSP, we can use DSP to solve
\eqref{e-tractable-classification} and \eqref{e-tractable-regression}
with the same convenience and ease. All that is required is importing DSP via
\verb|from dsp import *|, and 
replacing lines 8-16 above with the code block below.

\begin{lstlisting}[language=mypython,firstnumber=8]
# Creating local variables for uncertainty set
x1, x2 = LocalVariable(d), LocalVariable(d)
	
# Create bilinear form of theta and x
g1, g2 = saddle_inner(theta, x1), saddle_inner(-theta, x2)
		
# Construct uncertainty set containing x (filled in by user)
local_constraints1 = []
local_constraints2 = []

# Take suprema over x
G1 = saddle_max(g1, local_constraints1)
G2 = saddle_max(g2, local_constraints2)
\end{lstlisting}
For the user's convenience, in the Appendix we present
a helper Python function that automatically converts
problems of the form \eqref{e-main} to problems of the form
\eqref{e-tractable-classification} and \eqref{e-tractable-regression}.

\section{Example}
\label{s-experiments}

We consider the problem of predicting nightly Airbnb rental prices in
London, from different features such as coordinates, distance from city
center, and neighborhood restaurant quality index.
We will consider a simulated hypothetical case where we do not have full acess
to the rentals' location. 
We will use robust regression to handle the uncertain
location features.
We can then use uncertainty sets for the unknown
locations, allowing us to illustrate the ease of
specifying RERM problems with our framework. 
This example is artificial, but does use real original data.
We do not advocate using robust regression in particular for this
problem; replacing each unknown location with a center of the uncertainty set
performs nearly as well as the best robust regression method, and is much
simpler. The code to reproduce this example is available at
\begin{quote}
    \url{https://github.com/cvxgrp/rerm_code}.
\end{quote}

\paragraph{Data.} We begin with a curated dataset from London \cite{gyodi2021airbnb},
and remove rentals with prices exceeding 1000~Euros and those located more than 
7~km from the city center, resulting in a dataset of 3400 rows and 20
columns.  We then remove categorical features and randomly sub-sample
to obtain a training set with 1000 data points and test data set with 500
data points.
The training feature matrix is $X \in \reals^{1000 \times 9}$, with rows $x_i^T$.
The first two columns of $X$ correspond to the longitude and latitude
respectively.
Our baseline predictor of rental price is a simple linear ordinary least squares (OLS)
regression based on all~9 features. Its test RMS error is 138 euros.

\paragraph{Hidden location features.}

\begin{figure}
	\centering
	\includegraphics[width=10cm]{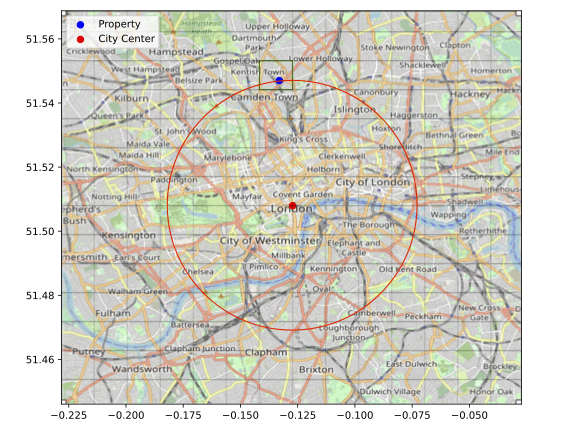}
	\caption{A visualization of $D_i$ and $S_i$ for a particular
	rental $i$. The large
	disk around the red dot corresponds to $D_i$.
	The square containing the blue dot corresponds to $S_i$.
	The overlap of the square and the disk corresponds to $D_i \cap S_i$.}
	\label{fig-london-map}
\end{figure}

To illustrate our method, we imagine a case where rental owners have elected to not
release the exact longitude and latitude of their properties.
(While not identical to this example, Airbnb does in fact mask rental locations.)
We grid London into 1~km by 1~km squares, and for each rental only give
the square $S_i$ it is located in.
This generates an uncertainty set for data point $i$ given by
\[
	\mathcal{X}_i^S = \{ x \in \reals^9 \mid x_{1:2} \in S_i,~ x_{3:9} =
X^i_{3:9} \}.
\]
We also know the distance of each rental from the city center $c\in\reals^2$,
denoted by $d_i$. Using this, we can consider a more refined 
uncertainty set
\[
	\mathcal{X}_i^{S\cap D} = \mathcal{X}_i^S \cap D_i,
\]
where $D_i=\{x\in\reals^9 \mid\|x_{1:2}-c\|_2\leq d_i\}$. See
figure~\ref{fig-london-map} for a visualization of these uncertainty sets.

We solve \eqref{e-tractable-regression} with
squared loss $f = \vert \cdot \vert^2$ and the two choices of the uncertainty sets
described above. These choices correspond to using
square or disk-intersected-with-square uncertainty sets for the missing
coordinates. Note that the square uncertainty set combined with the quadratic loss is a 
special case where we can derive an analytical form for the worst-case loss.
However, the analytical form is lost once we intersect the square with the disk. 

\paragraph{Comparing the methods.}

\begin{figure}
\centering
\includegraphics[width=8cm]{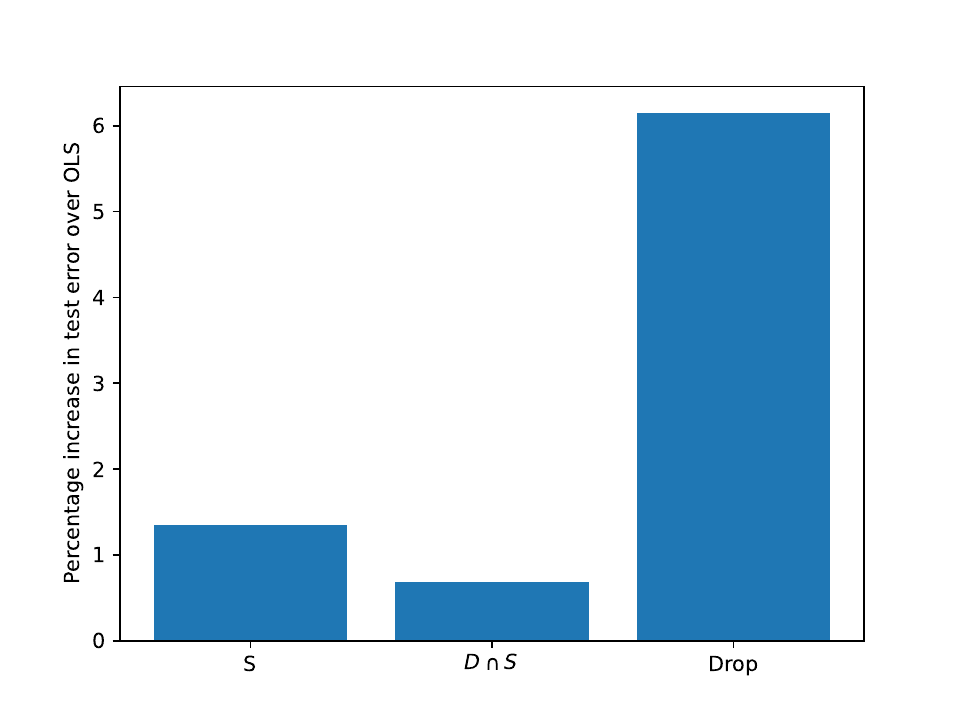}
\caption{Excess test error in the Airbnb price
prediction
experiment.}
\label{fig-real-exp1}
\end{figure}

We depict the performance of the two RERM predictors, as well as a predictor that
completely ignores coordinate information, in figure~\ref{fig-real-exp1}.
Our performance metric is the mean squared error on the test set,
in excess of the baseline OLS predictor trained on $X$ without any missing entries.

We observe that the dropping scheme, denoted as Drop in
figure~\ref{fig-real-exp1},
performs the worst. 
The robust predictors that use square uncertainty sets (denoted
as $S$) and the intersected uncertainty sets (denoted as
$S \cap D$) outperform
the others. The robust predictor that uses uncertainty sets $S \cap D$ outperforms
the robust predictor that uses only $S$.
Indeed, its performance is nearly as good as
the OLS predictor baseline that has access to all the columns
of $X$. These intersected uncertainty sets are
complex and do not admit the sort of convenient reformulation
afforded by using square or disk uncertainty sets. Yet, our framework allows us
to handle these uncertainty sets conveniently, and hence obtain
less conservative predictors.

\section*{Acknowledgements}
Stephen Boyd was partially supported by ACCESS (AI Chip Center for Emerging
Smart Systems), sponsored by InnoHK funding, Hong Kong SAR, and by Office of
Naval Research grant N00014-22-1-2121. 
This material is based upon work supported by the National Science
Foundation Graduate Research Fellowship Program under
Grant No. DGE1745016. Any opinions, findings, and conclusions or
recommendations expressed in this material are those of the authors
and do not necessarily reflect the views of the National Science Foundation.
This work was also supported by the AI Institute for Resilient Agriculture
(NSF-USDA award number 2021-67021-35329) 
\printbibliography

\clearpage
\appendix
\section*{Appendix}
In this section, we present a helper function that automatically converts
problems of the form \eqref{e-main} to problems of the form
\eqref{e-tractable-classification} and \eqref{e-tractable-regression}. This
helper function requires the user to pass as inputs the CVXPY variables,
constraints and loss function that define \eqref{e-main}, and is
presented below. The code is available at
\begin{quote}
    \url{https://github.com/cvxgrp/rerm_code}.
\end{quote}

\begin{lstlisting}[language=mypython, deletekeywords={abs}]
from typing import Callable, List
import cvxpy as cp
import numpy as np
from cvxpy.transforms.suppfunc import SuppFunc

def form_rerm(
    f: Callable,
    y: np.ndarray,
    theta: cp.Variable,
    theta_constraints: List[cp.Constraint],
    xs: List[cp.Variable],
    x_constraints: List[List[cp.Constraint]],
    mode: str
):
    """
    Args:
        f: a convex function
        y: a vector of length n
        theta: a CVXPY variable
        theta_constraints: a list of constraints on theta
        xs: a list of n scalar CVXPY variables
        x_constraints: a list of n lists of constraints on xs
        mode: "non_increasing" or "non_decreasing_sym_abs"

    Returns:
        A CVXPY problem instance of the robust ERM
        problem. 
    """
    n = len(xs)
    assert theta.ndim <= 1
    assert n == len(x_constraints) == len(y)

    z = cp.Variable(n)
    obj = 0.0
    constraints = theta_constraints
    
    for i in range(n):
        obj += f(z[i])
        G = SuppFunc(xs[i], x_constraints[i])
        if mode == "non_increasing":
            constraints += [
                G(-theta) + y[i] <= -z[i],
            ]
        elif mode == "non_decreasing_sym_abs":
            constraints += [
                G(theta) - y[i] <= z[i],
                G(-theta) + y[i] <= z[i],
            ]
        else:
            raise NotImplementedError

    prob = cp.Problem(cp.Minimize(obj), constraints)
    return prob
\end{lstlisting}
We emphasize that the loss
function $f$ provided by the user is not verified to have the
claimed curvature properties specified in the mode argument. It is
impossible to verify this in general, so the user must be careful to
provide a loss function that has the correct curvature
properties. We also mention
that the user is not limited to pass in a loss function $f$ that is a CVXPY
atom. Instead, the user has the flexibility to create a loss function
that is a composition of several CVXPY atoms, as follows.

\begin{lstlisting}[language=mypython]
def f(x: cp.Variable):
    """
    An example non-decreasing convex function of the magnitude of x.
    This is how a user can specify a arbitrary convex function.
    """
    return cp.square(cp.power(cp.abs(x), 1.5))
\end{lstlisting}

\end{document}

%% file: defs.tex
\newcommand{\reals}{{\mbox{\bf R}}}

\newcommand{\ie}{{\it i.e.}}

\newcommand{\BEAS}{\begin{eqnarray*}}
\newcommand{\EEAS}{\end{eqnarray*}}
\newcommand{\BEA}{\begin{eqnarray}}
\newcommand{\EEA}{\end{eqnarray}}
\newcommand{\BEQ}{\begin{equation}}
\newcommand{\EEQ}{\end{equation}}
\newcommand{\BIT}{\begin{itemize}}
\newcommand{\EIT}{\end{itemize}}

\newcommand{\todo}[1]{\textcolor{blue}{\textbf{EL: #1}}}
\newcommand{\dm}[1]{\textcolor{red}{\textbf{DM: #1}}}

\newcounter{algorithmctr}
\renewcommand{\thealgorithmctr}{\arabic{algorithmctr}}
   {\mbox{}\\*[\parskip]\begin{minipage}{\linewidth}%
       \refstepcounter{algorithmctr}\begin{list}{}{%
       \setlength{\rightmargin}{0\linewidth}%
       \setlength{\leftmargin}{.05\linewidth}}%
       \rmfamily\small
       \item[]{\setlength{\parskip}{0ex}\hrulefill\par%
        \nopagebreak{\bfseries\textsf{Algorithm \thealgorithmctr~}}}}%
   {{\setlength{\parskip}{-1ex}\nopagebreak\par\hrulefill\\*[2ex]\par}%
   \end{list}\end{minipage}}